\documentclass[twoside]{article}
\usepackage{graphicx}
\usepackage{amsmath}
\usepackage{theorem}
\usepackage{amsfonts}
\usepackage{amssymb}
\input psfig.sty         
 \theoremstyle{change}
\newtheorem{theorem}{Theorem}[section]
\newtheorem{proposition}[theorem]{Proposition}
{ \theorembodyfont{\upshape}

\newtheorem{example}[theorem]{Example}
}

{\ \rule{0.5em}{0.5em}\smallskip}

\begin{document}
\newcommand{\minisecless}[1]{\par\smallskip\noindent{\bf #1}\,}
\newcommand{\minisec}[1]{\par\bigskip\noindent{\bf #1}\,}
\renewcommand{\div}{\operatorname{div}}
\newcommand{\im}{\operatorname{im}}
\newcommand{\vol}{\operatorname{vol}}
\newcommand{\Sym}{\operatorname{Sym}}
\newcommand{\tr}{\operatorname{tr}}
\newcommand{\eps}{\varepsilon}
\newcommand{\pa}{\partial}
\newcommand{\RP}{\mathbb{RP}}
\newcommand{\CP}{\mathbb{CP}}
\newcommand{\CC}{\mathbb{C}}
\newcommand{\RR}{\mathbb{R}}
\newcommand{\ZZ}{\mathbb{Z}}
\newcommand{\MMM}{{\cal M}}
\newcommand{\WWW}{{\cal W}}
\newcommand{\DDD}{{\cal D}}
\newcommand{\LLL}{{\cal L}}
\newcommand{\FFF}{{\cal F}}
\newcommand{\U}{U}

\raggedbottom

\title{Gaussian densities and stability for some Ricci solitons}
\author{Huai-Dong Cao\footnote{
With partial support from NSF Grant DMS-0206847.}\\Lehigh
University \and Richard Hamilton\\Columbia University \and Tom
Ilmanen\footnote{ With partial support from Schweizerische
Nationalfonds Grant 21-66743.01.}\\ETH-Z\"{u}rich and Columbia
University}
\date{April 7, 2004}
\maketitle

\bigskip
Perelman \cite{Perelman} has discovered a remarkable variational
structure for the Ricci flow: it can be viewed as the gradient
flow of the entropy functional $\lambda$.  There are also two
monotonicity formulas of shrinking or localizing type: the
shrinking entropy $\nu$, and the reduced volume.  Either of these
can be seen as the analogue of Huisken's monotonicity formula for
mean curvature flow \cite{Huisken-shrink}. In various settings,
they can be used to show that centered rescalings converge
subsequentially to shrinking solitons, which function as idealized
models for singularity formation.

In this note, we exhibit the second variation of the $\lambda$ and
$\nu$ functionals, and investigate the linear stability of
examples. We also define the ``central density'' of a shrinking
Ricci soliton ({\it shrinker}) and compute its value for certain
examples in dimension 4. Using these tools, one can sometimes
predict or limit the formation of singularities in the Ricci flow.
In particular, we show that certain Einstein manifolds are
unstable for the Ricci flow in the sense that generic
perturbations acquire higher entropy and thus can never return
near the original metric. A detailed version of the calculations
summarized in this announcement will follow in
\cite{Cao-Hamilton-Ilmanen}.

In \S\ref{SecVarSteady}, we investigate the stability of
Perelman's $\lambda$-functional.  Its critical points are steady
solitons (Ricci flat in the compact case). We compute the second
variation $\DDD^2\lambda$; the corresponding Jacobi field operator
$L$ is a degenerate negative elliptic integro-differential
operator. In fact, $L$ equals half the Lichnerowicz Laplacian
$\Delta_L$ on divergence-free symmetric tensors, and zero on Lie
derivatives.  This fact and further investigations of the second
variation have been reported by Perelman \cite{Perelman-talk}. We
call a steady soliton {\it linearly stable} if $L\le0$, otherwise
{\it linearly unstable}. If $g$ is linearly unstable, then $g$ can
be perturbed so that $\lambda(g)>0$, which will destabilize it
utterly: it will decay into a cacophony of shrinkers and disappear
in finite time.  One observes that $\lambda(g)\le0$ for any metric
on the torus $T^n$, $n\le7$; in fact this is equivalent to the
positive mass theorem. By Guenther, Isenberg and Knopf
\cite{Guenther-Isenberg-Knopf} every $K3$ surface is linearly
stable; more generally, by Dai, Wang and Wei \cite{Dai-Wang-Wei}
any manifold with a parallel spinor is linearly stable. Other
cases are open.

In \S\ref{SecVarShrink}, we investigate the stability of the
$\nu$-functional, whose critical points are shrinkers.  The Jacobi
field operator $N$ of $\nu$ is like $L$ but with lower order
terms. We call a shrinker {\it linearly stable} if $N\le0$. Again,
$N$ is closely related to the Lichnerowicz Laplacian.  We observe
that $\CP^N$ (with the standard metric) is linearly stable, but
all other compact complex surfaces with $c_1>0$ are linearly
unstable.  Using results of Gasqui and Goldschmidt
\cite{Gasqui-Goldschmidt-QQQ,Gasqui-Goldschmidt-QQQQ} the complex
hyperquadric $Q^3$ (a Hermitian symmetric space) is linearly
unstable. This implies that $Q^3$ is irremediably unstable in the
sense that a generic (non-K\"ahler!) perturbation of $Q^3$ will
never approach the original geometry of $Q^3$ at any scale or
time. On the other hand, the hyperquadric $Q^4$ is linearly
stable. Other cases are open.

A notion of central density (or gaussian density) of a shrinker
can be defined from either of Perelman's monotonicity formulas; we
call these notions $\Theta$ and $\nu$. On a shrinker, the two
definitions are equivalent via $\Theta=e^\nu$. For a general
solution, $e^\nu$ is a lower bound for the central density of any
shrinker that arises later as a singularity model, which restricts
the shrinkers that may occur later.  This is presented in
\S\ref{SetupSec}.

The central density of certain standard 4-dimensional examples are
exhibited in a table in \S\ref{TableSec}.

We are grateful to Hugh Bray, Robert Bryant, Hubert Goldschmidt,
Dan Knopf, and John Morgan for illuminating conversations.

\section{Second Variation of the Entropy $\lambda$}
\label{SecVarSteady}

The second variation of the Einstein functional is positive in the
conformal direction but negative in all other directions.  The
glory of Perelman's entropy is that there is a preliminary
minimization over scalar functions that absorbs nearly all the
positive directions: the Jacobi field operator is a linear
integro-differential operator with nonpositive symbol.\footnote{It
is convenient that this happens without changing the metric via a
conformal change.  So the scalar variations do not affect the
background geometry. Instead they satisfy a {\it linear} PDE.}

Fix a compact manifold $(M,g)$. Define \cite{Perelman}
$$
\FFF(g,f):=\int e^{-f}(|Df|^2+R).
$$
Define the {\it entropy}
$$
\lambda(g):=\inf\{\FFF(g,f):f\in C^\infty_c(M), \int e^{-f}=1\}
$$
The infimum is achieved by a function $f$ solving
$$
-2\Delta f+|Df|^2-R=\lambda(g).
$$
Now consider variations $g(s)=g+sh$. Following Perelman, the first
variation $\DDD_g\lambda(h)$ of $\lambda$ is given by
\begin{align*}
\left.\frac{d}{ds}\right|_{s=0}\lambda(g(s))= \int
e^{-f}(-Rc-D^2f):h,
\end{align*}
where $f$ is the minimizer. A stationary point satisfies
$$
Rc+D^2f=0,
$$
which implies that $g$ is a {\it (gradient) steady soliton}, that
is, the Ricci flow with initial condition $g$ satisfies
$$
g(t)=\phi_t^*(g)
$$
where $\phi_t$ is a family of diffeomorphism generated by the
gradient vector field $Df$. In fact, any compact steady is Ricci
flat with $f=0$, $\lambda=0$.

Note by diffeomorphism invariance of $\lambda$ that
$\DDD_g\lambda$ vanishes on any Lie derivative $h=L_xg$.  From
this, by inserting $h=-2(Rc+D^2f)$ one recovers Perelman's
wonderful result that $\lambda(g(t))$ is nondecreasing on any
Ricci a Ricci flow, and is constant if and only if $g(t)$ is a
steady soliton.

We prove the following. Write $Rm(h,h):=R_{ijkl}h_{ik}h_{jl}$,
$\div\omega:=D_i\omega_i$, $(\div h)_i:=D_jh_{ji}$,
$(\div^*\omega)_{ij}=-(D_i\omega_j+D_j\omega_i)/2
=-(1/2)L_{\omega^\#}g_{ij}$.

\begin{theorem} The second variation
$\DDD^2_g\lambda(h,h)$ of $\lambda$ on a compact Ricci flat
manifold is given by
\begin{align*}
\left.\frac{d^2}{ds^2}\right|_{s=0}\lambda(g(s)) &=\int
-\frac{1}{2}|Dh|^2+|\div h|^2-\frac{1}{2}|Dv_h|^2
+Rm(h,h)\\
&=\int Lh:h,
\end{align*}
where
$$
Lh:=\frac{1}{2}\Delta h+\div^*\div h
+\frac{1}{2}D^2v_{h}+Rm(h,\cdot),
$$
and $v_h$ satisfies
\begin{align*}
\Delta v_h=\div\div h.
\end{align*}
\end{theorem}

\noindent The symbol of $L$ in the direction $\xi\in T_x^*M$ is
$$
\sigma_\xi(h)=-\pi_{\xi^\perp}(h),
$$
where $\pi_{\xi^\perp}(h)$ restricts $h$ to the hyperplane
$\xi^\perp$. So the operator $L$ is degenerate negative elliptic,
and has a discrete spectrum with at most a finite-dimensional
space of positive eigenfunctions.

Decompose $C^\infty(\Sym^2(T^*M))$ as
$$
\ker\div\oplus\im\div^*.
$$
One verifies that $L$ vanishes on $\im\div^*$, that is, on Lie
derivatives. On $\ker\div$ one has
$$
L=\frac{1}{2}\Delta_L
$$
where
$$
\Delta_Lh:=\Delta h+2Rm(h,\cdot)-Rc\cdot h-h\cdot Rc
$$
is the Lichnerowicz Laplacian on symmetric 2-tensors.

We call a critical point $g$ of $\lambda$ {\it linearly stable} if
$L\le0$, and a {\it maximizer} if $\lambda(g_1)\le\lambda(g)$ for
all $g_1$. A compact Ricci flat metric is a maximizer if and only
if it admits no metric of positive scalar curvature. (This follows
from Schoen's solution of the Yamabe problem \cite{Schoen}.)
Evidently a maximizer is stable.  If $g$ is not stable, then a
slight perturbation will develop $\lambda>0$ and $R>0$ and (in
principle) disappear in finite time as positive manifolds do. A
good question is whether any Ricci-flat manifold is unstable.  We
call this the {\it positive mass problem for Ricci flat
manifolds.}

\begin{example}
$T^n$ admits no metric of positive scalar curvature by the
positive mass theorem, so $\lambda(g)\le0$ for all $g$ on $T^n$.
\end{example}

\begin{example}
A Calabi-Yau K3 surface and more generally, any manifold with a
parallel spinor has $\Delta_L\le0$
\cite{Guenther-Isenberg-Knopf,Dai-Wang-Wei}. So these manifolds
are linearly stable in the sense presented here.
\end{example}

\begin{example}
Let $g$ be compact and Ricci flat. Following \cite{Buzzanca,
Guenther-Isenberg-Knopf} we examine conformal variations. It is
convenient to replace $ug$ by
$$
h=Su:=(\Delta u)g-D^2u
$$
which differs from the conformal direction only by a Lie
derivative and is divergence free. We have
$$
\Delta_LSu=(S\Delta u)g,
$$
so $\Delta_L$ has the same eigenvalues as $\Delta$. In particular,
$N\le0$ in the conformal direction. This contrasts with the
Einstein functional.
\end{example}

\section{Second Variation of the Shrinker Entropy $\nu$}
\label{SecVarShrink}

Fix a complete manifold $(M,g)$. Define
$$
\WWW(g,f,\tau):=\frac{1}{(4\pi\tau)^{n/2}}\int e^{-f}\left[
\tau(|Df|^2+R)+f-n\right]\,dV.
$$
Define the {\it shrinker entropy} by
$$
\nu(g):=\inf\{\WWW(g,f,\tau):f\in C^\infty_c(M), \tau>0,
\frac{1}{(4\pi\tau)^{n/2}}\int e^{-f}=1\}
$$
Assume that $M$ is compact or is asymptotic at infinity to a
metric cone over a smooth, compact Riemannian manifold.  One
checks that $\nu(g)$ is realized by a pair $(f,\tau)$ that solve
the equations
\begin{align*}
\tau(-2\Delta f+|Df|^2-R)-f+n+\nu=0,\qquad
\frac{1}{(4\pi\tau)^{n/2}}\int fe^{-f}=\frac{n}{2}+\nu,
\end{align*}
and $f$ grows quadratically.

Consider variations $g(s)=g+sh$ where $h$ is smooth of compact
support.  Following Perelman, one calculates the first variation
$\DDD_g\nu(h)$ to be
\begin{align*}
\left.\frac{d}{ds}\right|_{s=0}\nu(g(s))=
\frac{1}{(4\pi\tau)^{n/2}}\int e^{-f}(\tau(-Rc-D^2f)+g/2):h.
\end{align*}
A stationary point of $\nu$ satisfies
\begin{align}
D^2f+Rc-\frac{g}{2\tau}=0
\end{align}
which says that $g$ is a (gradient) shrinker, that is, its Ricci
flow $g(t)$ has the form
\begin{align*}
g(t):=(T-t)\psi_t^*(g)),\qquad t<T,
\end{align*}
where $\psi_t$ are the diffeomorphisms generated by $-Df$, and
$\tau=T-t$.

As before, $\DDD_g\nu$ vanishes on Lie derivatives. By scale
invariance it vanishes on multiplies of the metric. Inserting
$h=-2(Rc+D^2f-g/2\tau)$, one recovers Perelman's brilliant formula
that finds that $\nu(g(t))$ is monotone on a Ricci flow, and
constant if and only if $g(t)$ is a gradient shrinker.

A positive Einstein manifold is a shrinker with $f\equiv n/2$,
normalized by $Rc=g/2\tau$. We compute:

\begin{theorem} Let $(M,g)$ be a positive Einstein manifold.
The second variation $\DDD^2_g\nu(h,h)$ is given by
\begin{align*}
\left.\frac{d^2}{ds^2}\right|_{s=0}\nu(g(s))
&=\frac{\tau}{\vol(g)}\int-\frac{1}{2}|Dh|^2+|\div&h|^2-\frac{1}{2}|Dv_h|^2
+Rm(h,h)+\frac{v_h^2}{4\tau}\\
&&-\frac{1}{2n}\left(\frac{1}{\vol(g)}\int \tr_gh\right)^2,\\
&=\frac{\tau}{\vol(g)}\int Nh:h,&\\
\end{align*}
where
$$
Nh:=\frac{1}{2}\Delta h+\div^*\div h+
\frac{1}{2}D^2v_h+Rm(h,\cdot)-\frac{g}{2n\tau\vol(g)}\int\tr_gh.
$$
and $v_h$ is the unique solution of
\begin{align*}
\Delta v_h+\frac{v_h}{2\tau}=\div\div h,\qquad\int v_h=0.
\end{align*}
\end{theorem}

\noindent There is a strictly more complicated formula in the case
of non-Einstein shrinkers.

As in the previous case, $N$ is degenerate negative elliptic and
vanishes on $\im\div^*$.  Write
$$
\ker\div=(\ker\div)_0\oplus\RR g
$$
where $(\ker\div)_0$ is defined by $\int\tr_gh=0$. Then on
$(\ker\div)_0$ we have
$$
N=\frac{1}{2}\left(\Delta_L-\frac{1}{\tau}\right)
$$
where $\Delta_L$ is the Lichnerowicz Laplacian. So the linear
stability of a shrinker comes down to the (divergence free)
eigenvalues of the Lichnerowicz Laplacian. Let us write $\mu_L$
for the maximum eigenvalue of $\Delta_L$ on symmetric 2-tensors
and $\mu_N$ for the maximum eigenvalue of $N$ on $(\ker\div)_0$,

\begin{example} The round sphere is {\it geometrically stable}
(i.e. nearby metrics are attracted to it up to scale and gauge) by
the results of Hamilton
\cite{Hamilton-three,Hamilton-four,Hamilton-two} and Huisken
\cite{Huisken-ricci}. In particular it is linearly stable:
$\mu_N=-2/(n-1)\tau<0$.
\end{example}

\begin{example}
\label{CPn-example} For $\CP^N$, the maximum eigenvalue of
$\Delta_L$ on $(\ker\div)_0$ is $\mu_L=1/\tau$ by work of
Goldschmidt \cite{Gasqui-Goldschmidt-talk}, so $\CP^N$ is
neutrally linearly stable, i.e. the maximum eigenvalue of $N$ on
$(\ker\div)_0$ is $\mu_N=0$.
\end{example}

Any product of two nonflat shrinkers $N_1^{n_1}\times N_2^{n_2}$
is linearly unstable, with $\mu_N=1/2\tau$. The destabilizing
direction $h=g_1/n_1-g_2/n_2$ corresponds to a growing discrepancy
in the size of the factors.

More generally, any compact K\"ahler shrinker with $\dim
H^{1,1}(M)\ge2$ is linearly unstable. Again, this can be seen
directly: a small perturbation into a non-canonical K\"ahler class
will move in a straight line nearly toward the vertex of the
K\"ahler cone, hence away from the canonical class (in a scale
invariant sense).  If $M$ is K\"ahler-Einstein, we compute $\mu_N$
as follows. Let $\sigma$ be a harmonic 2-form and $h$ be the
corresponding metric perturbation; then $\Delta_Lh=0$, and if
$\sigma$ is chosen perpendicular to the K\"ahler form, then as
above we obtain $\mu_N=1/2\tau$.

A complete list of compact complex surfaces with $c_1>0$ is
$\CP^2$, $\CP^1\times\CP^1$, and $\CP^2\#k(-\CP^2)$,
$k=1,\ldots8$. Each of these has a unique K\"ahler shrinker metric
(K\"ahler-Einstein unless $k=1,2$). By the above, all are linearly
unstable except $\CP^2$.

Let $Q^N$ denote the complex hyperquadric in $\CP^{N+1}$ defined
by
$$
\sum_{i=0}^{N+1}z_i^2=0,
$$
a Hermitian symmetric space of compact type, hence a positive
K\"ahler-Einstein manifold. Then $Q^2$ is isometric to
$\CP^1\times\CP^1$, the simplest example of the above instability
phenomenon.

\begin{example}
Consider $Q^3$. It has $\dim H^{1,1}(Q^3)=1$, so the above
discussion does not apply. But the maximum eigenvalue of
$\Delta_L$ on $(\ker\div)_0$ is $\mu_L=-2/3\tau$ by work of Gasqui
and Goldschmidt \cite{Gasqui-Goldschmidt-QQQ} (or see
\cite{Gasqui-Goldschmidt}). The proximate cause is a
representation that appears in the sections of the symmetric
tensors but not in scalars or vectors. Therefore, $Q^3$ is
linearly unstable with
$$
\mu_N=\frac{1}{6\tau}.
$$
Since $\nu$ increases along the Ricci flow, this implies that a
generic small perturbation of the Einstein metric $\bar g$ will
grow and $g(t)$ will never return near $\bar g$ at any scale or
time. We say that $\bar g$ is {\it geometrically unstable}.  Now
imagine that we start with a random metric on $Q^3$ and propose to
use the Ricci flow to find an Einstein metric or other canonical
geometry for $Q^3$.  Assuming there are no other critical points,
we find that the flow combusts in singularities of more elementary
type and the topology of the underlying manifold simplifies
drastically, unless it happens to get hung up at the Einstein
metric $\nu$. So the Ricci flow has fundamentally more complicated
behavior than in dimension three, as one expects. Further
exploration of this example will appear in
\cite{Bryant-Goldschmidt-Ilmanen-Morgan}.
\end{example}

\begin{example}
Let $Q^4$ be the 4-dimensional hyperquadric. The maximum
eigenvalue of $\Delta_L$ on symmetric tensors is $\mu_L=-1/\tau$
by work of Gasqui and Goldschmidt \cite{Gasqui-Goldschmidt-QQQQ}
(or see \cite{Gasqui-Goldschmidt}). So $Q^4$ is neutrally linearly
stable: $\mu_N=0$.
\end{example}

Let $g$ be a positive Einstein metric, and let us examine
conformal variations. As before, without loss replace $ug$ by the
divergence-free variation
$$
h=Su:=(\Delta u)g-D^2u+\frac{ug}{2\tau}.
$$
As before, $\Delta_LSu=(S\Delta u)g$. Thus $\Delta_L$ has the same
eigenvalues as $\Delta_{\text{fns}}|(\ker S)^\perp$. But $\ker S$
is empty except on round $S^n$, which is linearly stable.  Note
that $\mu_{\text{fns}}\le -n/2(n-1)\tau$ with equality only on
round $S^n$. So we have:

\begin{proposition}
A positive Einstein metric is linearly unstable for conformal
variations if and only if the maximum eigenvalue of $\Delta$ on
functions satisfies
$$
-\frac{1}{\tau}< \mu_{\text{fns}}< -\frac{n}{2(n-1)\tau}.
$$
\end{proposition}

\noindent We do not know whether this inequality can ever be
satisfied on a positive Einstein manifold.

\section{The Central Density of a Shrinker}
\label{SetupSec}

Our aim in this section is to define the central density of a
gradient shrinker.  First we define a suitable class of gradient
shrinkers, then we review the two Perelman monotonicity formulas
of shrinking type and apply them by taking the center point to be
the parabolic vertex of the shrinker. Our principal result is that
the two notions of density coincide.

A gradient shrinker solves
\begin{align*}
\pa g/\pa t=-2Rc,\qquad g(t):=-t\psi_t^*(g(-1)),\qquad t<0,
\end{align*}
where $\psi_t$ are the diffeomorphisms generated by the gradient
of a function $F(x,t)$. Differentiating the above expression
yields
\begin{align}
\label{Shrinker} D^2F+Rc-\frac{g}{2\tau}=0,
\end{align}
where $\tau=-t$. Normalizing $F$ by adding a time-dependent
constant, we obtain
\begin{align}
\label{ShrinkerMove} \frac{\pa F}{\pa\tau}+|DF|^2=0.
\end{align}
Differentiating (\ref{Shrinker}), taking the trace two ways, and
applying Bianchi II and commutation rules yields
$D(|DF|^2+R-F/\tau)=0$. Adding a further global constant to $F$
leads to the classical {\it auxiliary equation}
\begin{align}
\label{ShrinkerExtra} |DF|^2+R-\frac{F}{\tau}=0.
\end{align}
Equations (\ref{Shrinker})-(\ref{ShrinkerExtra}) are the
fundamental equations for a gradient shrinker. Combining
(\ref{Shrinker}) and (\ref{ShrinkerMove}) yields the backward heat
equation
\begin{equation}
\label{ShrinkerParaF} \frac{\pa F}{\pa\tau}=\Delta
F-|DF|^2+R-\frac{n}{2\tau}.
\end{equation}

In order to prove our results we need some analytic hypotheses on
the metric of $M$. We assume that $M$ is complete, connected, and
$\kappa$-noncollapsed at all scales. We also assume that the
curvature decays quadratically as $x\to\infty$. (This is
satisfied, for example, by the blowdown shrinker $L(N,-1)$
\cite{Feldman-Ilmanen-Knopf}.) However, many of our results hold
under the weaker hypothesis of bounded curvature.  Under the
quadratic decay hypothesis, $g(t)$ converges in the
Gromov-Hausdorff sense as $t\nearrow0$ to a metric cone $C$ which
is smooth except at the vertex, which we call $0$.  The
convergence is smooth except on a compact set, which falls into
the vertex, which we call $0$. For a proof and further details,
see \cite{Ilmanen-shrinker}.

\bigskip
We now wish to define a gaussian density centered at the parabolic
vertex $(y,s)=(0,0)$ of the spacetime
$\MMM:=(M\times(-\infty,0))\cup (C\times\{0\})$. We may do this in
two ways: via the shrinking entropy or the reduced volume, both
due to Perelman \cite{Perelman}.

The reduced volume generalizes Bishop volume monotonicity to the
spacetime setting. For a smooth point $(y,s)$ in any Ricci flow,
define the {\it reduced distance} $\ell=\ell_{y,s}$ by
\begin{align}
\ell(x,t)&:=\frac{1}{2\sqrt{\tau}}\inf_\gamma\int_0^\tau\sqrt{\sigma}
\left(\left|\frac{d\gamma}{du}\right|^2+R\right)\,d\sigma,\qquad
t<s,\quad x\in M,
\end{align}
where the infimum is taken over all paths $(\gamma(u),u)$, $t\le
u\le s$ that connect $(x,t)$ to $(y,s)$. The {\it reduced volume
centered at $(y,s)$} is defined by
$$
\theta_{y,s}(t):=\frac{1}{(4\pi\tau)^{n/2}}\int_Me^{-\ell(x,t)}\,dV_t(x).
$$
Perelman wonderfully shows that $\theta_{y,s}(t)$ is increasing in
$t$ and is constant precisely on a gradient shrinker. Now define
$\ell_{0,0}$ by passing smooth points $(y_i,s_i)$ to $(0,0)$. We
have:

\begin{proposition}
$\ell_{0,0}$ is well-defined and is locally Lipschitz on $\MMM$.
For $t<0$, $\theta_{0,0}(t)$ is well defined, constant, and
contained in (0,1].
\end{proposition}

\noindent This constant value we call the {\it central density of
(M,g(t))} and denote
$$
\Theta(M)=\Theta(M,g(\cdot)):=\theta_{0,0}(t),\qquad t<0.
$$

Next we turn to the shrinking entropy $\nu$. Let $(M,g(t))$ be a
smooth Ricci flow existing up to $t=s$ and set $\tau:=s-t$. Let
$f$ solve the heat equation
$$
\frac{\pa f}{\pa t}=\Delta f-|Df|^2+R-\frac{n}{2\tau},
$$
that came up for the soliton potential of a shrinker
(\ref{ShrinkerParaF}). Define $u$ by
$$
u:=\frac{e^{-f}}{(4\pi\tau)^{n/2}}.
$$
Remarkably, $u$ solves the adjoint heat equation
\begin{align}
\label{ShrinkerParaU} \frac{\pa u}{\pa t}=\Delta u-Ru.
\end{align}
This leads to the conservation law
$$
\frac{1}{(4\pi\tau)^{n/2}}\int e^{-f}=\int u=1\qquad \text{for
}t<s.
$$
Perelman has shown \cite{Perelman}
$$
\frac{\pa}{\pa t}\WWW(s-t,f(t),g(t)) =\frac{1}{(4\pi\tau)^{n/2}}
\int2e^{-f}\left|D^2f+Rc-\frac{g}{2\tau}\right|^2\,dV_t\ge0.
$$
and the right hand side vanishes precisely when $g(t)$ is a
gradient shrinker and $f$ is its soliton potential. (This shows,
as mentioned above, that $\nu$ increases in general and is
constant on a shrinker.)

If $u$ emerges from a dirac source at a smooth point $(y,s)$, we
write $u=u_{y,s}$, $f=f_{y,s}$ and define the {\it shrinker
entropy centered at $(y,s)$} by
$$
\phi_{y,s}(t):=\WWW(s-t,f_{y,s}(t),g(t))
$$
Passing smooth points $(y_i,s_i)$ to $(0,0)$, we prove:

\begin{proposition}
$u_{0,0}$ is well-defined, smooth, and positive on
$M\times(-\infty,0)$ and solves equation (\ref{ShrinkerParaU}). It
satisfies
\begin{equation}
\label{PerelNormal} \int_Mu_{0,0}\equiv
\frac{1}{(4\pi\tau)^{n/2}}\int e^{-f_{0,0}}=1,\qquad t<0.
\end{equation}
Also, $\phi_{0,0}(t)$ is well-defined and lies in $(-\infty,0]$
for all $t<0$. In fact, it is constant, with
$$
\phi_{0,0}(t)=\nu(M),\qquad t<0.
$$
\end{proposition}

Since $\phi_{0,0}(t)$ is constant, $f_{0,0}$ is a soliton
potential and so
\begin{equation}
\label{PlusC} f_{0,0}=F+C
\end{equation}
for some constant $C$ depending only on $\MMM$.

We now wish to relate $\Theta(M)$ and $\nu(M)$ via $F$.  In the
process we determine the value of $C$, and sharpen on a shrinker
the general Perelman relation $\ell_{0,0}\le f_{0,0}$
\cite{Perelman}.

We begin with $\Theta(M)$. Using the symmetry of the shrinker and
a simple comparison argument, one checks:

\begin{proposition}
The integral curves of $F$ are minimizing $\LLL$-geodesics
emanating from $(0,0)$.
\end{proposition}

Then using the homothetic time-symmetry of the shrinker, one
obtains after a straightforward computation:
$$
\ell_{0,0}=\tau(|DF|^2+R)=F=f_{0,0}-C,
$$
and thus by the definition and (\ref{PerelNormal}), one gets:
$$
\Theta(M)=e^C.
$$
Next, we evaluate $\nu$. Compute
\begin{align*}
\nu(M)&=\phi_{0,0}(t)\\
&=\frac{1}{(4\pi\tau)^{n/2}}\int e^{-f_{0,0}}\left[
\tau(|Df_{0,0}|^2+R)+f_{0,0}-n\right]\,dV_t\\
&=\frac{1}{(4\pi\tau)^{n/2}}\int e^{-f_{0,0}}\left[
\tau(\Delta f_{0,0}+R)+f_{0,0}-n\right]\,dV_t\\
&=\frac{1}{(4\pi\tau)^{n/2}}\int e^{-f_{0,0}}\left[
f_{0,0}-n/2\right]\,dV_t
\end{align*}
by integrating by parts and the trace of (\ref{Shrinker}).  On the
other hand, by (\ref{ShrinkerExtra}) and (\ref{PlusC}), the first
integral expression also equals
\begin{align*}
\nu(M)&=\frac{1}{(4\pi\tau)^{n/2}}\int e^{-f_{0,0}}\left[
f_{0,0}-C+f_{0,0}-n\right]\,dV_t
\end{align*}
We conclude that $\nu(M)=2\nu(M)-C$, so $\nu(M)=C$. We summarize
these results in a theorem.

\begin{theorem}
\label{Log} On a shrinker satisfying the above assumptions, we
have
$$
\ell_{0,0}=F=f_{0,0}-\nu(M)
$$
and
$$
\Theta(M)=e^{\nu(M)}.
$$
\end{theorem}

Details will appear in \cite{Cao-Hamilton-Ilmanen}.

\section{Table of Values}
\label{TableSec}

In this section we calculate $\Theta$ for some standard shrinkers.
The computations are simplified by several observations. Normalize
positive Einstein manifolds by $Rc=g/2\tau$, $\tau=1/2(n-1)$, so
that $S^n$ has radius $1$.

(1) For any shrinker $M$, $\Theta(M)\le1$ with equality if and
only if $M=\RR^n$.

(2) Let $M$ be a Ricci flat cone with $g=dr^2+r^2g_\Sigma$ where
$\Sigma$ is positive Einstein. Then $M$ is a shrinker (with
interior singularity), and
$$
\Theta(M)=\frac{\vol(\Sigma)}{\vol(S^n)}.
$$

(3) If $M$ is a positive Einstein manifold, then (for any $\tau$)
$$
\Theta(M)=\left(\frac{1}{4\pi\tau e}\right)^{n/2}\vol_\tau(M)\le
\Theta(S^n),
$$
with equality if and only if $M=S^n$.

(4) $\Theta(S^n)=\left(\dfrac{n-1}{2\pi e}\right)^{n/2}\vol(S^n)$.
By way of comparison, note that for mean curvature flow,
$\Theta_{MCF}(S^n)=\left(\dfrac{n}{2\pi
e}\right)^{n/2}\vol(S^n)$.\footnote{Following an observation of
White, we note (tantalizingly) that the respective limits as
$n\to\infty$ are $\sqrt{2/e}$ and $\sqrt{2}$.}

(5) $\Theta(\CP^N)=\left(\dfrac{N+1}{\pi
e}\right)^N\dfrac{\vol(S^{2N+1})}{2\pi}$.

(6) The positive K\"ahler-Einstein manifold $M=\CP^2\#k(-\CP^2)$,
$k=0,3,\ldots,8$,  has $\Theta(M)=(9-k)/2e^2$.

(7) $\Theta(M\times N)=\Theta(M)\Theta(N)$.

We say that one shrinker {\it decays} to another if there is a
small perturbation of the first whose Ricci flow develops a
singularity modelled on the second.  Because the $\nu$-invariant
is monotone during the flow, decay can only occur from a shrinker
of lower density to one of higher density. This creates a ``decay
lowerarchy''. (It should be a partial order.)

We have computed the following density values in dimension 4. Note
that the conclusion of Theorem \ref{Log} holds for all our
examples, though not all are smooth enough to satisfy the
hypotheses.

\bigskip
\begin{tabular}
[c]{|lllr|} \hline
&&&\\
\textbf{Shrinker} & \textbf{Type} & $\Theta$ & $\Theta$ \\
\hline
&&&\\
$\RR^4$ & flat & 1 & 1.000\\
&&&\\
$S^4$ & positive Einstein & ${6/e^2}$ & .812\\    
&&&\\
$S^3\times\RR$ & product & ${2\left(\pi/e^3\right)^{1/2}}$ & .791\\
&&&\\
$S^2\times\RR^2$ & product & ${2/ e}$ & $.736$\\  
&&&\\
$L(2,-1)$ & blowdown shrinker \cite{Feldman-Ilmanen-Knopf}&
${e^{\sqrt{2}-2}(1+\sqrt{2})/2}$ & .672\\   
&&&\\
$\CP^2$ & positive Einstein & ${9/2e^2}$ & .609\\  
&&&\\
$S^2\times S^2$ & product & ${4/ e^2}$ & .541 \\ 
&&&\\
$\CP^2\#(-\CP^2)$ & Koiso metric \cite{Koiso,Cao} & ${3.826/e^2}$ & .518\\
&&&\\$\CP^2\#(-\CP^2)$ & Page metric \cite{Page} & ${3.821/e^2}$ & .517\\
&&&\\
$C(\RP^3)$ & Ricci flat cone & ${1/2}$ & .500 \\
&&&\\
$C(\RP^2)\times\RR$ & product & ${1/ 2}$ & .500 \\
&&&\\
$\RP^4$ & positive Einstein & ${3/ e^2}$ & .406\\   
&&&\\
$\CP^2\#3(-\CP^2)$ & positive Einstein & ${3/e^2}$ & .406\\  
&&&\\
$\RP^3\times\RR$ & product & $\left(\pi/ e^3\right)^{1/2}$ & .396 \\
&&&\\
$\RP^2\times\RR^2$ & product & ${1/ e}$ & .368 \\ 
&&&\\
$\CP^2\#4(-\CP^2)$ & positive Einstein & ${5/2e^2}$ & .338\\  
&&&\\
$C(S^3/\ZZ_3)$ & Ricci flat cone & ${1/ 3}$ & .333 \\
&&&\\
$\CP^2\#5(-\CP^2)$ & positive Einstein & ${2/e^2}$ & .271\\  
&&&\\
$\CP^2\#6(-\CP^2)$ & positive Einstein & ${3/2e^2}$ & .203\\  
&&&\\
$\CP^2\#7(-\CP^2)$ & positive Einstein & ${1/e^2}$ & .135\\  
&&&\\
$\CP^2\#8(-\CP^2)$ & positive Einstein & ${1/2e^2}$ & .068\\  
&&&\\
\hline
\end{tabular}
\bigskip

All manifolds in the table are created from Einstein manifolds
except for $L(2,-1)$ and the Koiso metric. The computations for
these metrics will be detailed in \cite{Cao-Hamilton-Ilmanen}. The
volume of the Page metric is computed in \cite{Page}.

The {\it blowdown shrinker} $L(n,-1)$ is a K\"ahler shrinker
defined on the total space of the tautological holomorphic line
bundle $o_{N-1}(-1)$ over $\CP^{N-1}$, that is, on $\CC^N$ blown
up at $z=0$. The metric of $L(N,-1)$ is $U(N)$ invariant,
complete, and conelike at infinity, satisfying quadratic decay for
the curvature. As $t\nearrow0$, the exceptional divisor
$\CP^{N-1}$ shrinks to a point and elsewhere the metric converges
smoothly to a cone metric on $\CC^N\setminus\{0\}$ whose metric
completion has a vertex at $0$. For positive time, the flow can
continue by a smooth, $U(N)$-invariant K\"ahler expander on
$\CC^N$. See [FIK].

The Koiso metric \cite{Koiso,Cao} and the Page metric \cite{Page}
are both $\U(2)$-invariant metrics on $\CP^2\#(-\CP^2)$.  The
former, but not the latter, is K\"ahler.  The remarks following
Example \ref{CPn-example} show that the Koiso metric has one
direction of instability (in a K\"ahler direction). On the other
hand, the Page metric may well decay to the Koiso metric. By the
discussion in \cite{Feldman-Ilmanen-Knopf}, this leads us to
conjecture that either metric decays to $\CP^2$ via a $\CP^1$
pinches off.

\end{document}